\documentclass[final,5p,times,twocolumn]{elsarticle}

\usepackage{amssymb}
\usepackage{amsmath}

\newcommand{\mb}[1]{\mbox{\boldmath $#1$}}
\newcommand{\mbs}[1]{{\mbox{\boldmath \scriptsize{$#1$}}}}
\usepackage{amsthm}
\newtheorem{theorem}{Theorem}
\newtheorem{definition}[theorem]{Definition}
\newtheorem{lemma}[theorem]{Lemma}
\newtheorem{example}{Example}
\newtheorem{conj}{Conjecture}
\usepackage{tikz}
\usepackage{caption,subcaption}
\usetikzlibrary{shapes,arrows,positioning, intersections}

\tikzset{vertex/.style = {shape=circle,draw,minimum size=1.5em}}
\tikzstyle{line} = [draw, -latex']

\begin{document}

\begin{frontmatter}
\title{Pairwise independent correlation gap}

\author[au1]{Arjun Ramachandra} 
\address[au1]{{Decision Sciences Area, Indian Institute of Management Bangalore},
            {Bannerghatta Road},
      {Bangalore},
          {560076},
         {India}}
      \ead{arjun.ramachandra@iimb.ac.in}
            \author[au2]{Karthik Natarajan\fnref{l2}} 
\fntext[l2]{Corresponding author}
\address[au2]{{Engineering Systems and Design, Singapore University of Technology and Design},
           {8 Somapah Road},
            {Singapore},
			{487372},
	{Singapore}}
      \ead{karthik\_natarajan@sutd.edu.sg}
            
\begin{abstract}
In this paper, we introduce the notion of a ``pairwise independent correlation gap'' for set functions with random elements. The pairwise independent correlation gap is defined as the ratio of the maximum expected value of a set function with arbitrary dependence among the elements with fixed marginal probabilities to the maximum expected value with pairwise independent elements with the same marginal probabilities. We show that for any nonnegative monotone submodular set function defined on $n$ elements, this ratio is upper bounded by $4/3$ in the following two cases: (a) $n = 3$ for all marginal probabilities and (b) all $n$ for small marginal probabilities (and similarly large marginal probabilities). This differs from the bound on the ``correlation gap'' which holds with mutual independence and showcases the fundamental difference between pairwise independence and mutual independence. We discuss the implication of the results with two examples and end the paper with a conjecture.
\end{abstract}




\end{frontmatter}

\section{Introduction}
\label{submission}
Submodular set functions play an important role in optimization and machine learning \cite{Lovasz,bachsub}. An important notion which characterizes the behavior of such functions with random elements is the ``correlation gap'', which was introduced in \cite{Agrawal2012}, building on the work of \cite{Calinescu2}. The correlation gap is defined as the ratio of the maximum expected value of a set function (not necessarily submodular) with arbitrary dependence among the elements with fixed marginal probabilities to the expected value of the function with mutually independent elements with the same marginal probabilities. A seminal result derived in \cite{Calinescu2,Agrawal2012} showed that for any nonnegative monotone submodular set function, the correlation gap is upper bounded by $e/(e-1) \approx 1.582$. Correlation gap has found applicability in distributionally robust optimization \cite{Agrawal2012}, content resolution schemes \cite{rico}, mechanism design \cite{yan} and combinatorial prophet inequalities \cite{rubinstein}. In this paper, we introduce the notion of the ``pairwise independent correlation gap'' which is defined as the ratio of the maximum expected value of a set function with arbitrary dependence among the elements with fixed marginal probabilities to the maximum expected value of the function with pairwise independent elements with the same marginal probabilities. While mutual independence implies pairwise independence, the reverse implication does not hold beyond two random variables. Interestingly, we show that for nonnegative monotone submodular set functions, the pairwise independent correlation gap is upper bounded by $4/3 \approx 1.334$ in two cases, showcasing the difference between pairwise independence and mutual independence.

\subsection{Preliminaries}
Let $N$ denote a ground set with $|N| = n$ elements and $f: 2^{N} \rightarrow \mathbb{R}_+$ be a nonnegative set function. The function is submodular if $f(S) + f(T) \geq f(S \cap T) +  f(S \cup T)$ for all $S, T \subseteq N$ or equivalently $f(S \cup \{i\}) - f(S) \geq f(T \cup \{i\}) -  f(T)$ for all $S \subseteq T$, $i \in N \backslash T.$ The function is monotone (non-decreasing) if $ f(S) \leq f(T)$ for all $S \subseteq T$. Examples of nonnegative monotone submodular set functions include the rank function of matroids, the coverage function and the entropy function \cite{bachsub}. The multilinear extension of the set function is given by $F: [0,1]^n \rightarrow \mathbb{R}_+$ where $F(\mb{x})$ is the expected value of the function where each element $i \in N$ is selected with probability $x_i$ independently:
\begin{equation}
\begin{array}{rlllll}\label{eq:ind}
\displaystyle F(\mb{x})  = \sum_{S \subseteq N}\prod_{i \in S}x_i \prod_{i \in N \backslash S}{(1-x_i)}f(S).
\end{array}
\end{equation}
The concave closure of the set function is given by $f^+: [0,1]^n \rightarrow \mathbb{R}_+$ where $f^{+}(\mb{x})$ is the maximum expected value of the function over all joint distributions of the elements where each $i \in N$ is selected with probability $x_i$:
\begin{equation}
\begin{array}{rlllll}\label{eq:univprimal}
\displaystyle f^{+}(\mb{x}) = \max \Big{\{}\sum_{S \subseteq N} \theta(S)f(S) \ \Big{|} \ \sum_{S \subseteq N: S \ni i} \theta(S) = x_i,  \forall i \in N, \\
\displaystyle \sum_{S \subseteq N} \theta(S) = 1, \theta(S) \geq 0,  \forall S \subseteq N\Big{\}}.
\end{array}
\end{equation}
Both the multilinear extension and the concave closure define continuous extensions of $f$ such that  $0 \leq F(\mb{x}) \leq f^{+}(\mb{x})$ for all $\mb{x} \in [0,1]^n$. The correlation gap is defined as the ratio ${f^{+}(\mb{x})}/{F(\mb{x})}$, where we assume ${0}/{0}=1$.
\begin{theorem}
\label{thm:calinescushipra} \cite{Calinescu2,Agrawal2012}
For any $n$, any nonnegative monotone submodular function $f: 2^{N} \rightarrow \mathbb{R}_+$ and any $\mb{x} \in [0,1]^n$, ${f^{+}(\mb{x})}/{F(\mb{x})} \leq e/(e-1)$. This bound is sharp and attained when $f(S) =\min (|S|,1)$ and $\mb{x} =(1/n,\ldots,1/n)$ as $n \rightarrow \infty$.
\end{theorem}
For fixed $n$, the analysis in \cite{rico} sharpens the bound to $1/(1-(1-1/n)^n)$ which asymptotically converges to $e/(e-1)$. Specifically for $n = 2$, the bound on the correlation gap is $4/3$ and for $n = 3$, the bound is $27/19$.

\section{Set Functions with Pairwise Independence}
Given Bernoulli random variables $(\tilde{x}_1,\ldots,\tilde{x}_n)$, mutual independence is equivalent to the condition $P(\tilde{x}_i = x_i; i \in N) = \prod_{i=1}^{n}P(\tilde{x}_i = x_i)$ for all $\mb{x} \in \{0,1\}^n$. Pairwise independence is a weaker notion where only pairs of random variables are independent, namely $P(\tilde{x}_i = x_i,\tilde{x}_j = x_j) = P(\tilde{x}_i = x_i)P(\tilde{x}_j = x_j)$ for all $(x_i,x_j)\in \{0,1\}^2$ and $i < j, \; i,j \in N$. One of the motivations for studying pairwise independent random variables is that the joint distribution can have a low cardinality support (quadratic in $n$), in contrast to mutually independent random variables (exponential in $n$). The low cardinality of these distributions have important ramifications in efficiently derandomizing algorithms for NP-hard optimization problems \cite{luby}. Pairwise independence has also proved to be useful in modeling data with zero or weak correlations but more complex higher-order dependencies. Such behavior has been observed experimentally in cortical neurons \cite{schenid}. Recently there has been an interest in understanding the behavior of stochastic optimization problems by relaxing the assumption of mutual independence to pairwise independence \cite{ramachandra2021pairwise,cara,dughmi,gupta}. Our results are related to this line of work.

 \subsection{A Pairwise Independent Continuous Extension}
Here we consider pairwise independent random elements to model limited independence.
\begin{definition}
\label{def:pairwiseind}
The maximum expected value of a set function over all pairwise independent distributions where each element $i \in N$ is selected with probability $x_i$ is given by:
\begin{equation}
\begin{array}{rlllll}\label{eq:bivprimal}
\displaystyle f^{++}(\mb{x}) = \max & \displaystyle \sum_{S \subseteq N} \theta(S)f(S)  \\
  \mbox{s.t} & \displaystyle \sum_{S \subseteq N: S \ni i, j} \theta(S) = x_ix_j, & \forall i < j, i,j \in N\\
     & \displaystyle \sum_{S \subseteq N: S \ni i} \theta(S) = x_i, & \forall i \in N  \\
  & \displaystyle \sum_{S \subseteq N} \theta(S) = 1 &  \\
  & \theta(S) \geq 0, & \forall S \subseteq N.
\end{array}
\end{equation}
We refer to $f^{++}(\mb{x})$ as the pairwise independent extension of a set function and define the pairwise independent correlation gap as the ratio ${f^{+}(\mb{x})}/{f^{++}(\mb{x})}$.
\end{definition}
Clearly, we have $F(\mb{x}) \leq f^{++}(\mb{x}) \leq f^{+}(\mb{x})$ where $F(\mb{x}) = f^{++}(\mb{x})$ when $n = 2$. In the special case when $f(S) =\min (|S|,1)$, \cite{ramachandra2021pairwise} recently showed that $f^{++}(\mb{x}) = \min(\sum_ix_i-\max_i x_i\; (\sum_ix_i-\max_ix_i),1)$. In this case, $f^{+}(\mb{x})= \min(\sum_i x_i,1)$ and $F(\mb{x}) = 1-\prod_i(1-x_i)$ giving rise to an upper bound of $4/3$ for the pairwise independent correlation gap and $1/(1-(1-1/n)^n)$ for the correlation gap. A natural question is whether such bounds generalize to other nonnegative monotone submodular functions? We answer this question in part in this paper by showing that the $4/3$ upper bound holds for all nonnegative monotone submodular functions in two cases: (a) $n = 3$ with general marginal probabilities (the first nontrivial case where pairwise and mutual independence are different notions) and (b) general $n$ with small (and similarly large) marginal probabilities. Our proofs are based on constructions of pairwise independent distributions. Other such constructions can be found in \cite{geisser1962,karloffmansour1994,derriennic2000}. Our results essentially show that the extremal pairwise independent distribution is provably closer to the extremal distribution with arbitrary dependence in comparison to the mutually independent distribution.

\section{Upper Bound for $n = 3$}\label{sec:smallngeneralx}
We assume throughout that $f$ is monotone submodular. Furthermore, let $f(\emptyset) = 0$ and $f(N) = 1$. This is without loss of generality, since translating and scaling the function preserves monotonicity and submodularity. Let ${\cal F}_n$ denote the set of all such functions defined over $n$ elements:
\begin{equation*}
\begin{array}{rllll}
{\cal F}_n = \{f  \ | \ f \mbox{ is} \mbox{ monotone submodular}, f(\emptyset) = 0, f(N) = 1\}.
\end{array}
\end{equation*}
Observe that ${\cal F}_n$ is a polytope defined over $2^n$ variables where the variables correspond to the set function values. Unfortunately, computing $f^{+}(\mb{x})$ and likewise $f^{++}(\mb{x})$ is challenging for functions in the set ${\cal F}_n$ (in fact NP-hard in the worst-case; see \cite{Agrawal2012}). In particular, the probability distributions which attain the values $f^{+}(\mb{x})$ and $f^{++}(\mb{x})$ are not oblivious to the function $f \in {\cal F}_n$. This is in contrast to $F(\mb{x})$ which is evaluated for the mutually independent distribution. We illustrate this in Tables \ref{tab:nonob} and \ref{tab:nonob1} using two functions in ${\cal F}_n$ for $n = 3$ with $(x_1,x_2,x_3) = (1/2,1/2,1/2)$.
\begin{table}[!htbp]
\caption{Optimal distribution of $f^{+}(\mb{x})$ and $f^{++}(\mb{x})$ for given $f$.}
\label{tab:nonob}
\begin{center}
\begin{tabular}{ccccc}
$S$ & $f(S)$ & $\theta^{+}(S)$ & $\theta^{++}(S)$  & $\theta^{I}(S)$ \\ \hline
$\emptyset$    & 0 & 0 & $1/4$ &$1/8$ \\
$\{1\}$ & ${1}/{3}$ & $1/2$ & $0$ &$1/8$ \\
$\{2\}$ & ${1}/{2}$ & 0 &  $0$ &$1/8$ \\
$\{3\}$ & ${1}/{2}$ & 0 & $0$ &$1/8$\\
$\{1,2\}$ & ${5}/{6}$ &     0&  $1/4$ &$1/8$  \\
$\{1,3\}$ & ${5}/{6}$    & 0 & $1/4$ &$1/8$ \\
$\{2,3\}$  & $1$ &     $1/2$ &   $1/4$&$1/8$ \\
$\{1,2,3\}$  & $1$  & 0 & $0$  &$1/8$ \\ \hline
\end{tabular}
\end{center}
\end{table}
\begin{table}[!htbp]
\caption{Optimal distribution of $f^{+}(\mb{x})$ and $f^{++}(\mb{x})$ for given $f$.}
\label{tab:nonob1}
\begin{center}
\begin{tabular}{ccccc}
$S$ &  $f(S)$ & $\theta^{+}(S)$ & $\theta^{++}(S)$ & $\theta^{I}(S)$ \\ \hline
$\emptyset$    & 0 & 0 & 0 &$1/8$\\
$\{1\}$ &  $1/3$ & 0 & $1/4$&$1/8$\\
$\{2\}$ & $1/2$ & 0 & $1/4$&$1/8$\\
$\{3\}$ &  $1/2$ & $1/2$ & $1/4$&$1/8$\\
$\{1,2\}$ &  $3/4$ & $1/2$ &0   &$1/8$  \\
$\{1,3\}$ &  $3/4$ & 0 & 0&$1/8$\\
$\{2,3\}$  & $3/4$ & 0 & 0 &$1/8$\\
$\{1,2,3\}$  & 1 & 0  & $1/4$&$1/8$\\ \hline
\end{tabular}
\end{center}
\end{table}
Here $\theta^{+}$ and $\theta^{++}$ are the optimal distributions that attain $f^{+}(\mb{x})$ and $f^{++}(\mb{x})$ which is found by solving linear programs and $\theta^{I}$ is the mutually independent distribution. Even in this example, we see that analyzing the ratio ${f^{+}(\mb{x})}/{f^{++}(\mb{x})}$ is nontrivial since the optimal distributions in the numerator and denominator change with $f$ unlike ${f^{+}(\mb{x})}/{F(\mb{x})}$ where only the optimal distribution in the numerator changes. In the first case, $f^{+}(\mb{x}) = 2/3$, $f^{++}(\mb{x}) = 2/3$ and $F(\mb{x}) = 5/8$ while in the second case, $f^{+}(\mb{x}) = 5/8$, $f^{++}(\mb{x}) = 7/12$ and $F(\mb{x}) = 55/96$. Our technique of analysis of the ratio for $n = 3$ is based on partitioning the hypercube $[0,1]^3$ into smaller regions such that for each region, we obtain lower bounds on $f^{++}(\mb{x})$ (by generating feasible pairwise independent distributions of the linear program in \eqref{eq:bivprimal}) and upper bounds on $f^{+}(\mb{x})$ (by generating feasible solutions of the dual of the linear program in \eqref{eq:univprimal}). Using this, we prove the pairwise independent correlation gap is bounded from above by $4/3$ in two important cases. The following two lemmas will be useful in proving the result.
\begin{lemma}
\label{lem:extremepts}
The extremal functions in ${\cal F}_3$ are given by $\mathcal{E}\left({\cal F}_3\right)=\{E_1,E_2,E_3,E_4,E_5,E_6,E_7,E_8\}$ where
\begin{equation*}
\begin{array}{rlll}
E_1 = (0,1,0,0,1,1,0,1), E_2 = (0,0,1,0,1,0,1,1), \\
 E_3 = (0,0,0,1,0,1,1,1), E_4 = (0,1,1,0,1,1,1,1),\\
E_5 = (0,1,0,1,1,1,1,1), E_6 = (0,0,1,1,1,1,1,1), \\
 E_7 = (0,1,1,1,1,1,1,1), E_8 = (0,\frac{1}{2},\frac{1}{2},\frac{1}{2},1,1,1,1),
\end{array}
\end{equation*} with the values representing $(f(\{\emptyset\})$, $f(\{1\})$, $f(\{2\})$, $f(\{3\})$, $f(\{1,2\})$, $f(\{1,3\})$, $f(\{2,3\})$, $f(\{1,2,3\}))$.
In addition, the extreme points of ${\cal F}_3^{1} = {\cal F}_3 \cap \{f|f(\{1\})+f(\{2,3\}) \geq f(\{2\})+f(\{1,3\}),f(\{1\})+f(\{2,3\}) \geq f(\{3\})+f(\{1,2\})\}$ are given by $\mathcal{E}\left({\cal F}_3^{1}\right)=\{E_1,E_2,E_3,E_4,E_5,E_7,E_8\}$, the extreme points of ${\cal F}_3^{2} ={\cal F}_3 \cap \{f|f(\{2\})+f(\{1,3\}) \geq f(\{1\})+f(\{2,3\}),f(\{2\})+f(\{1,3\}) \geq f(\{3\})+f(\{1,2\}))\}$ are given by $\mathcal{E}\left({\cal F}_3^{2}\right)=\{E_1,E_2,E_3,E_4,E_6,E_7,E_8\}$ and the extreme points of ${\cal F}_3^{3} ={\cal F}_3 \cap \{f|f(\{3\})+f(\{1,2\}) \geq f(\{1\})+f(\{2,3\}),f(\{3\})+f(\{1,2\}) \geq f(\{2\})+f(\{1,3\})\}$ are given by $\mathcal{E}\left({\cal F}_3^{3}\right)=\{E_1,E_2,E_3,E_5,E_6,E_7,E_8\}$.
\end{lemma}
It is straightforward to verify that the extreme points of ${\cal F}_3$ are given by the lemma by enumerating linearly independent active constraints of the polytope and solving for them. The reader can also verify this using the open-source software \emph{polymake} \cite{polymake2017} which provides all the extreme points of small dimensional polytopes. A similar analysis can be carried out for ${\cal F}_3^{1}$ to ${\cal F}_3^{3}$ which are also polytopes.
\begin{lemma}
The following two inequalities $(I_1)$ and $(I_2)$ hold:
\begin{equation*}
\begin{array}{rlll}
(I_1) & \alpha+\beta-4\alpha\beta \geq 0,  \forall  [\alpha,\beta] \in [0,1]^2: \alpha+\beta \leq 1.\\
(I_2) & 4\alpha+4\beta-4\alpha\beta-3 \geq 0,  \forall [\alpha,\beta] \in [0,1]^2: \alpha+\beta \geq 1.
\end{array}
\end{equation*}
\end{lemma}
The inequalities hold since $\alpha+\beta-4\alpha\beta = (1-\alpha-\beta)(\alpha+\beta) + (\alpha-\beta)^2 $ and $4\alpha+4\beta-4\alpha\beta-3 =  (\alpha+\beta-1)(3-\alpha-\beta) + (\alpha-\beta)^2$.
This brings us to our first key theorem.
\begin{theorem}
\label{thm:n=3}
For $n= 3$, any nonnegative monotone submodular function and any $\mb{x} \in [0,1]^3$, $f^{+}(\mb{x})/f^{++}(\mb{x}) \leq  {4}/{3}$. This bound is sharp and attained when $f(S) = \min(|S|, 1)$ and $(x_1,x_2,x_3) = (0,1/2,1/2)$.
\end{theorem}
\begin{proof}
We outline the main steps of the proof below:
\begin{enumerate}
 \item
 Partition the unit hypercube $[0,1]^3$ into $5$ regions as follows:
 \begin{figure}[!htbp]
 \centering
 \begin{tikzpicture}[thick,scale=0.255, every node/.style={transform shape}]
      \node[scale=2.5][ draw,rectangle, text width=2cm,align=center] (fl1) {$\mb{x} \in [0,1]^3$};
      \node[scale=2.5][xshift=-4 cm,yshift=-2cm,draw, rectangle, text width=2.2cm,align=center] (fl2) { $x_1+x_2 \leq 1$};
        \node[scale=2.5][xshift=-6cm,yshift=-4cm,draw,rectangle, text width=2.5cm,align=center] (fl3){$x_1+x_2+x_3 \leq 1$};
         \node[scale=1.5][xshift=-10 cm,yshift=-7.75 cm,draw, vertex] (init){$\mb{R_1}$};
        \node[scale=2.5][xshift=0cm,yshift=-4cm,draw,rectangle, text width=2.5cm,align=center] (fl4){$x_1+x_2+x_3 > 1$};
\node[scale=2.5][xshift=-2cm,yshift=-6cm,draw,rectangle, text width=2.2cm,align=center] (fl5){$x_1+x_3 \leq 1$};
 \node[scale=2.5][xshift=2cm,yshift=-6cm,draw,rectangle, text width=2.2cm,align=center] (fl6){$x_1+x_3 > 1$};
  \node[scale=1.5][xshift=3.5 cm,yshift=-11.25 cm,draw, vertex] (init){$\mb{R_4}$};
 \node[scale=2.5][xshift=-4cm,yshift=-8cm,draw,rectangle, text width=2.2cm,align=center] (fl7){$x_2+x_3 \leq 1$};
  \node[scale=1.5][xshift=-6.75 cm,yshift=-14.5 cm,draw, vertex] (init){$\mb{R_2}$};
 \node[scale=2.5][xshift=0cm,yshift=-8cm,draw,rectangle, text width=2.2cm,align=center] (fl8){$x_2+x_3 > 1$};
  \node[scale=1.5][xshift=0 cm,yshift=-14.5 cm,draw, vertex] (init){$\mb{R_3}$};
 \node[scale=2.5][xshift=4cm,yshift=-2cm,draw,rectangle, text width=2.2cm,align=center] (fl9){$x_1+x_2 > 1$};
  \node[scale=1.5][xshift=6.75 cm,yshift=-4.5 cm,draw, vertex] (init){$\mb{R_5}$};
    \path[line] (fl1.south) --node[right]{}(fl2.north);
        \path[line] (fl2.south) -- (fl3.north);
        \path[line](fl2.south)  --  (fl4.north);
	    \path[line] (fl4.south) --  (fl5.north);
        \path[line] (fl4.south) --  (fl6.north);
        \path[line] (fl5.south) --  (fl7.north);
        \path[line] (fl5.south) --  (fl8.north);
        \path[line] (fl1.south) --  (fl9.north);
    \end{tikzpicture}
\caption{Partition tree of the probabilities in $[0,1]^3$ with $5$ regions represented by the leaves}
  \label{fig:treepartition}
    \end{figure}
(a) Sort the marginal probabilities in increasing value such that $0 \leq x_1 \leq x_2 \leq x_3 \leq 1$. In this case, we have $x_1+x_2 \leq x_1+x_3 \leq x_2+x_3 \leq x_1+x_2+x_3$. (b) Divide the hypercube into five regions $R_1$ to $R_5$ based on the relative position of the value $1$ among these partial sums (see Figure \ref{fig:treepartition}).
\item Consider the lower bound on $f^{++}(\mb{x})$ (denoted by $\underline{f}^{++}(\mb{x})$) provided by two feasible pairwise independent probability distributions for the linear program in \eqref{eq:bivprimal} (see Table \ref{tab:pwindfeas}).
\begin{table}[!htbp]
\caption{Pairwise independent distributions.}
\label{tab:pwindfeas}
\begin{center}
\begin{tabular}{ccc}
$S$ & $\theta_1(S)$ & $\theta_2(S)$\\ \hline
$\emptyset$    & $(1-x_3)(1-x_1-x_2)$ & $(1-x_2)(1-x_3)$ \\
$\{1\}$ & $x_1(1-x_3)$ & $0$\\
$\{2\}$    & $x_2(1-x_3)$ &  $(1-x_3)(x_2-x_1)$\\
$\{3\}$ & $x_3(1-x_1-x_2)+x_1x_2$ & $(1-x_2)(x_3-x_1)$\\
$\{1,2\}$      & $0$&  $x_1(1-x_3)$      \\
$\{1,3\}$     & $x_1(x_3-x_2)$ & $x_1(1-x_2)$    \\
$\{2,3\}$      & $x_2(x_3-x_1)$ &   $x_2x_3-x_1(x_2+x_3-1)$    \\
$\{1,2,3\}$    & $x_1x_2$ & $x_1(x_2+x_3-1)$  \\ \hline
Feasible &  $x_1+x_2\leq 1$ &  $x_2+x_3\geq 1$ 
\end{tabular}
\end{center}
\end{table}
The probability distribution $\theta_1$ is feasible if and only if $x_1+x_2 \leq 1$ and we use the lower bound generated by this distribution in regions $R_1$ to $R_4$. The probability distribution $\theta_2$ is feasible if and only if $x_2+x_3 \geq 1$. Since $x_1+x_2 \geq 1$ implies $x_2+x_3 \geq 1$, we use the lower bound generated by this distribution in region $R_5$.
\item
Consider the upper bound on $f^{+}(\mb{x})$ (denoted by $\overline{f}^+(\mb{x})$)
\begin{table}[!htbt]
\caption{Dual feasible solutions of \eqref{eq:univdual}}
\label{tab:feas-d}
\begin{center}
\begin{tabular}{cl}
Dual Feasible &  $\lambda_i$ values\\ \hline
  $D_1, f \in {\cal F}_3$ & $\lambda_0 = 0, \lambda_ 1= f(\{1\}),$ \\
  &  $\lambda_2 = f(\{2\}), \lambda_3 = f(\{3\})$\\
  $D_2, f \in {\cal F}_3^{1}$ & $\lambda_0 = f(\{2\})+f(\{3\})-f(\{2,3\}),$\\
  &  $\lambda _1= f(\{1\}-f(\{2\})-f(\{3\})+f(\{2,3\}),$\\
  &  $\lambda_2 = f(\{2,3\})-f(\{3\}),$\\
  &  $\lambda_3 = f(\{2,3\})-f(\{2\})$ \\
  $D_3 ,f \in {\cal F}_3^{2}$ & $\lambda_0 = f(\{1\})+f(\{3\})-f(\{1,3\}),$\\
  &   $\lambda_1 = f(\{1,3\})-f(\{3\}),$  \\
  &  $\lambda_2 = f(\{2\})-f(\{1\})-f(\{3\})+f(\{1,3\}),$\\
  &  $\lambda_3 =  f(\{1,3\})-f(\{1\})$\\
  $D_4, f \in {\cal F}_3^{3}$ & $\lambda_0 = f(\{1\})+f(\{2\})-f(\{1,2\}),$\\
  &  $\lambda_1 = f(\{1,2\})-f(\{2\}),$ \\
   & $\lambda_2 = f(\{1,2\})-f(\{1\}),$\\
   &   $\lambda_3 =  f(\{3\})-f(\{1\})-f(\{2\})+f(\{1,2\})$\\
$D_5,f \in {\cal F}_3^{1}$ & $\lambda_0 = f(\{1\})+f(\{3\})-f(\{1,3\}),$\\
&  $\lambda_1 = f(\{1,3\})-f(\{3\}),$\\
&  $\lambda_2 = f(\{2,3\})-f(\{3\}),$\\
&  $\lambda_3 =  f(\{1,3\})-f(\{1\})$ \\
 $D_6, f \in {\cal F}_3^{2}$  &  $\lambda_0 = f(\{2\})+f(\{3\})-f(\{2,3\}),$\\
 &  $\lambda_1 = f(\{1,3\})-f(\{3\}),$\\
 &  $\lambda_2 = f(\{2,3\})-f(\{3\}),$\\
&   $\lambda_3 =  f(\{2,3\})-f(\{2\})$ \\
 $D_7,f \in {\cal F}_3^{3}$ & $\lambda_0 = f(\{1,2\})+2f(\{3\})-f(\{2,3\})-f(\{1,3\})$\\
 &  $\lambda_1 = f(\{1,3\})-f(\{3\}),$\\
 &  $\lambda_2 = f(\{2,3\})-f(\{3\}),$\\
 &  $\lambda_3 = f(\{1,3\})+f(\{2,3\})-f(\{1,2\})-f(\{3\})$\\
 $D_8,f \in {\cal F}_3^3$& $\lambda_0 = f(\{2\})+f(\{3\})-f(\{2,3\}),$\\
 &  $\lambda_1 = f(\{1,2\})-f(\{2\}),$ \\ 
 &  $\lambda_2 = f(\{2,3\})-f(\{3\})$,\\
 &  $\lambda_3 = f(\{2,3\})-f(\{2\})$ \\
 $D_9,f \in {\cal F}_3$& $\lambda_0 = f(\{1,2\})+f(\{1,3\})+f(\{2,3\})-2$ \\ 
 &  $\lambda_1 = 1-f(\{2,3\}), \lambda_2= 1-f(\{1,3\}),$\\
 &  $\lambda_3 = 1-f(\{1,2\})$ \\ \hline
\end{tabular}
\end{center}
\end{table}
attained by a feasible solution of the dual linear program of \eqref{eq:univprimal}. The dual is given by:
\begin{equation}
\begin{array}{rllll}\label{eq:univdual}
\displaystyle f^{+}(\mb{x}) = \min \left\{\sum_{i\in N} \lambda_{i}x_i +\lambda_{0} \ \Big{|} \ \sum_{i \in S} \lambda_{i} +\lambda_{0} \geq f(S), \forall S \subseteq N\right\}.
\end{array}
\end{equation}
 Table \ref{tab:feas-d} displays nine dual solutions whose feasibility depends on $f \in {\cal F}_3$ (unlike the pairwise independent distributions in Table \ref{tab:pwindfeas} whose feasibility depends on the marginal probability vector $\mb{x}$ only). Here we use $f_i(S)$ to denote $f(S \cup \{i\})-f(S)$ While it is easy to verify $D_1, D_9$ are feasible for any  $f \in {\cal F}_3$, the remaining dual solutions are feasible for any $f$ in exactly one of ${\cal F}_3^{1},\; {\cal F}_3^{2},\;{\cal F}_3^{3}$, where ${\cal F}_3^{1} \cup {\cal F}_3^{2} \cup {\cal F}_3^{3}={\cal F}_3$ and ${\cal F}_3^{1} \cap {\cal F}_3^{2} \cap {\cal F}_3^{3}=\emptyset$.
\item We next choose appropriate pairwise independent primal feasible distributions  and dual feasible solutions for each region as outlined in Table \ref{tab:feas-pd}. In addition, we list the inequalities used in the proofs in each region.
 \begin{table}[!htbp]
\caption{Pairwise independent primal feasible distributions and dual feasible solutions of concave closure.}
\label{tab:feas-pd}
\begin{center}
\begin{tabular}{cccc}
Region & Primal Feasible & Dual Feasible & Inequalities \\ \hline
$R_1$    & $\theta_1$ & $D_1$ & $(I_1)$\\
$R_2$ & $\theta_1$ & $D_2,D_3,D_4$ &$(I_1),(I_2)$\\
$R_3$    & $\theta_1$ &  $D_2,D_6,D_8$ &$(I_1),(I_2)$\\
$R_4$ & $\theta_1$ & $D_5,D_6,D_7$ &$(I_1),(I_2)$\\
$R_5$      & $\theta_2$&  $D_9$    & $(I_2)$ \\ \hline
\end{tabular}
\end{center}
\end{table}
Specifically, our choice of dual feasible solutions for regions $R_2,R_3,R_4$ ensures that their feasible regions partition ${\cal F}_3$.
\item
Finally, we compute $\underline{f}^{++}(\mb{x})$ and $\overline{f}^+(\mb{x})$ using the chosen distributions and dual solutions from Table \ref{tab:feas-pd} for each region and use the fact that $4{f}^{++}(\mb{x}) - 3{f}^+(\mb{x}) \geq 4\underline{f}^{++}(\mb{x}) - 3\overline{f}^+(\mb{x})$. It is thus sufficient to prove that the inequality $4\underline{f}^{++}(\mb{x}) - 3\overline{f}^+(\mb{x}) \geq 0$ holds for each region and all $f \in {\cal F}_3$. Note that for a fixed $\mb{x} \in [0,1]^3$, $4\underline{f}^{++}(\mb{x}) - 3\overline{f}^+(\mb{x}) $ is linear in $f$ and we can thus minimize over $f \in {\cal F}_3$ for all $\mb{x} \in R_{j}, \; j \in \{1,2,3,4,5\}$.
We further note that for any region $R_j$, $4\underline{f}^{++}(\mb{x}) - 3\overline{f}^+(\mb{x})$ remains linear in $f$ as long as it's components $\underline{f}^{++}(\mb{x})$ and $\overline{f}^+(\mb{x})$ remain fixed for every $\mb{x} \in R_{j}$. Hence it is sufficient to verify the inequality for the extremal submodular functions of the polytopes provided n Lemma \ref{lem:extremepts}.
 \end{enumerate}
We illustrate the step for region $R_1$ where we need to verify that:
\begin{equation*}
\underset{f \in \mathcal{E}\left({\cal F}_3\right)}{\min} \;\left(4\underline{f}^{++}(\mb{x}) - 3\overline{f}^+(\mb{x})\right) \geq 0, \; \forall \mb{x} \in R_1.
\end{equation*}
For extreme point $E_1$, we have:
   \begin{equation*}
\begin{array}{lll}
4\underline{f}^{++}(\mb{x}) - 3\overline{f}^+(\mb{x}) \\
= 4 \sum_S \theta_1(S){E_1}(S) - 3\sum_i E_1(\{i\}) x_i\\
 = 4(x_1(1-x_3)+x_1(x_3-x_2)+x_1x_2)-3x_1\\
 = x_1\\
 \geq 0.
\end{array}
\end{equation*}
A similar argument works for extreme points $E_2$ and $E_3$. For extreme point $E_4$:
   \begin{equation*}
\begin{array}{lll}
4\underline{f}^{++}(\mb{x}) - 3\overline{f}^+(\mb{x})\\
 = 4 \sum_S \theta_1(S){E_4}(S) - 3\sum_i E_4(\{i\}) x_i\\
 = 4(x_1(1-x_3)+x_2(1-x_3)+x_1(x_3-x_2)+x_2(x_3-x_1)+x_1x_2)\\
 \ \ -3(x_1+x_2)\\
 = x_1+x_2-4x_1x_2\\
 \geq 0,
\end{array}
\end{equation*}
where the last inequality follows from $(I_1)$ with $\alpha=x_1$ and $\beta=x_2$. A similar argument works for extreme points $E_5$ and $E_6$. For extreme point $E_7$:   
\begin{equation*}
\begin{array}{lll}
 4\underline{f}^{++}(\mb{x}) - 3\overline{f}^+(\mb{x}) \\
  = 4 \sum_S \theta_1(S){E_7}(S) - 3\sum_i E_7(\{i\}) x_i\\
= 4(1-(1-x_3)(1-x_1-x_2))-3(x_1+x_2+x_3)\\
 = x_1+x_2+x_3-4x_1x_3-4x_2x_3\\
 \geq 0,
\end{array}
\end{equation*}
where the last inequality follows from $(I_1)$ with $\alpha=x_1+x_2$ and $\beta=x_3$. For extreme point $E_8$:
   \begin{equation*}
\begin{array}{lll}
4\underline{f}^{++}(\mb{x}) - 3\overline{f}^+(\mb{x})\\
 = 4 \sum_S \theta_1(S){E_8}(S) - 3\sum_i E_8(\{i\}) x_i\\
 = 4(x_1+x_2+x_3-x_1x_2)-3(x_1+x_2+x_3)\\
  =\left( x_1+x_2-4x_1x_2\right)+x_3\\
 \geq 0,
\end{array}
\end{equation*}
where the last inequality follows from $(I_1)$.
The proofs for the remaining four regions $R_2-R_5$ can be done similarly and we leave it to the reader to verify this\footnote{Details can be found at the online appendix}.
\end{proof}
The next lemma provides an immediate extension of the result to a simple sum of monotone submodular functions defined on disjoint ground sets of size at most three.

\begin{lemma}
\label{lem:usefullemma}
Let $N_1 ,N_2, \ldots, N_m$ be a partition of $N$ such that $N_i$ is of size at most $2$ or $3$ for all $i \in M$ where $M=\{1,2,\ldots,m\}$ and $f_i:2^{N_i} \rightarrow \mathbb{R}_{+}$ be monotone submodular set functions. Define $f = f_1+\ldots+f_m:2^N \rightarrow \mathbb{R}_{+} $. Then $f^{+}(\mb{x})/f^{++}(\mb{x}) \leq  {4}/{3}$ for all $\mb{x} \in [0,1]^n$.
\end{lemma}
\begin{proof}
Let $f^{+}_i(\mb{x})$ and $f^{++}_i(\mb{x})$ be the concave closure and pairwise independent extension of $f_i$ for each $i \in M$. Since the functions are defined over disjoint sets, we have $f^{+}(\mb{x})= \;\sum_{i=1}^m  f^{+}_i(\mb{x})$ and $f^{++}(\mb{x})= \;\sum_{i=1}^m  f^{++}_i(\mb{x})$.
Together this gives:
$$\frac{f^{+}(\mb{x})}{f^{++}(\mb{x})} = \frac{\sum_{i=1}^m \; f^{+}_i(\mb{x})}{\sum_{i=1}^m \; f^{++}_i(\mb{x})} \leq \frac{4}{3},$$ using $f^{+}_i(\mb{x})/f^{++}_i(\mb{x}) \leq 4/3$ for all $i \in M$ (where the case $n = 2$ follows from \cite{rico} and $n = 3$ from Theorem \ref{thm:n=3}).
\end{proof}

\section{Upper Bound for General $n$ with All Small or All Large Values in $\mb{x}$}
Our technique of analysis for general $n$ is based on the use of two pairwise independent distributions first proposed in \cite{ramachandra2021pairwise}. We show that these distributions suffice to prove a $4/3$ upper bound for all functions in ${\cal F}_n$ when the marginal probabilities are all small or all large.
\begin{lemma}
\label{lem:hw}
\cite{ramachandra2021pairwise} For any $\mb{x} \in [0,1]^n$, sort the values as $0 \leq x_1 \leq \ldots \leq x_n \leq 1 $.
\begin{enumerate}
\item Suppose $\sum_{i=1}^{n-1} x_i \leq 1$. Then there always exists a pairwise independent distribution with the marginal probability vector $\mb{x}$ of the form shown in Table \ref{tab:hunterworsley} (above). Specifically, the joint probabilities are given by: (a) $\theta(\emptyset) = (1- \sum_{i=1}^{n-1}x_i)(1-x_{n})$, (b) $\theta(\{i\}) = x_i(1-x_n)$ for all $i < n$, (c) $\theta(S) = 0$ for all other $S$ with $n \notin S$, and (d) $\theta(S) \geq 0 $ for all $S \ni n$ such that $\sum_{S: n \in S}\theta(S) = x_n$, $\sum_{S: i,n \in S}\theta(S) = x_ix_n$ for $i < n$, $\sum_{S: i,j,n \in S}\theta(S) = x_ix_j$ for all $i<j<n$.
\item  Suppose $\sum_{i=2}^{n} x_i \geq n-2$. Then there always exists a pairwise independent distribution with the marginal probability vector $\mb{x}$ of the form shown in Table \ref{tab:hunterworsley} (below). Specifically, the joint probabilities are given by: (a) $\theta(N) =x_{1}\left(\sum_{i=2}^{n}x_i -(n-2)\right)$, (b) $\theta(N\setminus \{i\}) = x_1(1-x_i)$ for all $i > 1$, (c) $\theta(S) = 0$ for all other $S$ with $1 \in S$, and (d) $\theta(S) \geq 0 $ for all $1 \notin S$ such that $\sum_{S: 1 \notin S}\theta(S) = 1-x_1$, $\sum_{S: 1,i \notin S}\theta(S) = (1-x_1)(1-x_i)$ for $i > 1$, $\sum_{S: 1,i,j \notin S}\theta(S) = (1-x_i)(1-x_j)$ for all $1<i<j$.
    \end{enumerate}
\end{lemma}

\begin{table}[!htbp]
\caption{Pairwise independent distributions with $\sum_{i=1}^{n-1} x_i \leq 1$ (above) and $\sum_{i=2}^{n} x_i \geq n-2$ (below).}
\label{tab:hunterworsley}
\begin{center}
\begin{tabular}{cc}
$S$ &  $\theta(S)$  \\ \hline
$\emptyset$ &  $\left( 1- \sum_{i=1}^{n-1}x_i \right)\left(1-x_{n} \right)$      \\
$\{1\}$ & $x_{1}(1-x_{n})$   \\
$\{2\}$ & $x_{2}(1-x_{n})$      \\
 \vdots & \vdots     \\
$\{n-1\}$ & $x_{n-1}(1-x_{n})$   \\
$\mbox{All other }S \mbox{ with } n \notin S$ & 0 \\
$\mbox{All }S \mbox{ with } n \in S$ & $\theta(S)$  \\ \hline
\end{tabular}
\begin{tabular}{cc}
$S$ &  $\theta(S)$ \\ \hline
$N$&  $x_{1}\left(\sum_{i=2}^{n}x_i -(n-2)\right)$     \\
$N\setminus \{2\}$ & $x_{1}(1-x_{2})$  \\
$N\setminus \{3\}$ & $x_{1}(1-x_{3})$     \\
 \vdots & \vdots     \\
$N\setminus \{n\}$ & $x_{1}(1-x_{n})$    \\
$\mbox{All other }S \mbox{ with } 1 \in S$ & 0\\
$\mbox{All }S \mbox{ with } 1 \notin S$ & $\theta(S)$ \\ \hline
\end{tabular}
\end{center}
\end{table}
This brings us to the main result of this section.
\begin{theorem}
\label{thm:generaln}
For any $n$, any nonnegative monotone submodular function and any $\mb{x} \in [0,1]^n$, we have $f^{+}(\mb{x})/f^{++}(\mb{x}) \leq  {4}/{3}$ in the following two cases:
\begin{enumerate}
\item
$\sum_{i=1}^{n-1} x_i \leq 1$ and $x_n \leq {1}/{4}$ where $0 \leq x_1 \leq \ldots \leq x_n \leq 1 $ (small probabilities)
\item
$\sum_{i=2}^{n} x_i \geq n-2$ and $x_1 \geq {3}/{4}$  where $0 \leq x_1 \leq \ldots \leq x_n \leq 1 $ (large probabilities).
\end{enumerate}
\end{theorem}
\begin{proof}
For small probabilities (1), to compute the lower bound on $f^{++}(\mb{x})$, we use the pairwise independent distribution provided in Lemma \ref{lem:hw}(1). This distribution is feasible when $\sum_{i=1}^{n-1}x_i \leq 1$ and was shown to be attain the pairwise independent extension in \cite{ramachandra2021pairwise} for the set function $f(S) =\min (|S|,1)$. Despite the non-optimality of this distribution for general $f \in {\cal F}_n$, we show that it suffices for our purposes. From the non-decreasing property of $f$, we have $ f(S) \geq f(\{n\})$ for all $S \ni n$. This implies that for the pairwise independent distribution in Lemma \ref{lem:hw}(1) we have:
$$\displaystyle {f}^{++}(\mb{x}) \geq  \sum_{i=1}^{n-1}x_i(1-x_{n}) f(\{i\})+x_{n}f(\{n\}).$$ For the upper bound, we use the dual feasible solution $\lambda_0 = 0$ and $\lambda_i= f(\{i\})$ for all $i \in N$. This is feasible for the dual linear program in (\ref{eq:univdual}) for all functions $f \in {\cal F}_n$ and gives:
$$\displaystyle {f}^+(\mb{x}) \leq  \sum_{i=1}^{n} x_if(\{i\}).$$
Together we get
     \begin{equation*}\label{eq:delta_smallprob_generaln}
\begin{array}{lll}
   \displaystyle 4{f}^{++}(\mb{x}) - 3{f}^+(\mb{x}) \\
     \geq \displaystyle 4 \sum_{i=1}^{n-1}x_i(1-x_{n}) f(\{i\})+4x_{n}f(\{n\})-3\sum_{i=1}^{n} x_if(\{i\})\\
     =   \displaystyle \sum_{i=1}^{n-1}x_i (1-4x_{n}) f(\{i\})+x_{n}f(\{n\})\\
     \geq  0 \;\;\mbox{ [since } x_n \leq 1/4].
    \end{array}
    \end{equation*}
Next, for the large probabilities (2), using the distribution in Lemma \ref{lem:hw}(2) we get:
 $${f}^{++}(\mb{x}) \geq \sum_{i=2}^{n}x_1(1-x_{i}) f(N\setminus \{i\}),$$
where we ignore the part of the distribution with $1 \notin S$ in Table \ref{tab:hunterworsley} (below) since $ f(S) \geq f(\emptyset)=0$ for all $S:1 \notin S$ (from the non-decreasing property of $f$).
  For the upper bound, we use the dual feasible solution $\lambda_0 =\sum_{i=1}^n f(N\setminus \{i\})-(n-1) $ and $\lambda_i= 1-f(N\setminus \{i\})$ for all $i \in N$. This is dual feasible for the linear program in (\ref{eq:univdual}) for all functions $f \in {\cal F}_n$ since
     \begin{equation*}
\begin{array}{llll}
  \sum_{i \in S} \lambda_{i} +\lambda_{0} \\
   =   \sum_{i \in S} \; \left[1-f(N\setminus \{i\})\right]+\left[\sum_{i=1}^n f(N\setminus \{i\})-(n-1)\right] \\
  = |S|-  \sum_{i \in S} f(N\setminus \{i\}) +\left[ \sum_{i \in S \cup S^c} f(N\setminus \{i\})-(n-1)\right] \\
  \mbox{ [where } S^c=N\setminus S] \\
  =\sum_{i \in S^c} f(N\setminus \{i\})+\left[|S|-(n-1)]\right]\\
  \geq [ f(S)+|S^c|-1]+[|S|-(n-1)]\\
  \mbox{ [using submodularity of f]}\\
   =  f(S) \;\;\forall S \subseteq N.
  \end{array}
     \end{equation*}
 We then get ${f}^+(\mb{x})
 \leq \sum_{i=1}^{n} (1-x_i)f(N\setminus \{i\})+\left[\sum_{i=1}^n x_i -(n-1)\right]$. Together we get
     \begin{equation*}
\begin{array}{lll}
 \displaystyle 4{f}^{++}(\mb{x})- 3{f}^+(\mb{x})  \\
  \geq    \displaystyle 4\left(\sum_{i=2}{n}x_1(1-x_{i}) f(N\setminus \{i\})+x_{1}\left(\sum_{i=2}{n}x_i -(n-2)\right)\right) \\
   \displaystyle -3\left( \displaystyle \sum_{i=1}{n} (1-x_i)f(N\setminus \{i\})+\sum_{i=1}{n} x_i -(n-1)\right)\\
        \geq     \displaystyle \displaystyle (4x_1-3)\left[ \sum_{i=2}{n} (1-x_i)f(N\setminus \{i\})+\left(\sum_{i=2}{n} x_i -(n-2)\right)\right] \\
        \mbox{    [since } 0 \leq f(N\setminus \{1\}) \leq 1]\\
     \geq  0\\
      \mbox{ [since }\sum_{i=2}{n}x_i \geq n-2 \mbox{ and } x_1 \geq 3/4].
    \end{array}
    \end{equation*}
\end{proof}

\section{Two Examples}\label{subsec:droapplication}
In this section, we discuss two examples to illustrate the usefulness of the new bounds.
\begin{example}
\textit{Distributionally robust $k$-sum optimization}: Given nonnegative weights $c_1,\ldots,c_n$ and a set ${\cal Y} \subseteq \{0,1\}^n$, the deterministic $k$-sum combinatorial optimization problem for a fixed integer $k \in \{1,\ldots,n\}$ is formulated as:
$$\min_{\mbs{y} \in {\cal Y} \subseteq \{0,1\}^n}  \max_{S \subseteq N:|S| \leq k}\sum_{i \in S}c_iy_i .$$
For $k = 1$, this reduces to a bottleneck combinatorial optimization problem. Suppose the weights $\tilde{c}_i$ are random where $\tilde{c}_i = c_i$ with probability $x_i$ and $\tilde{c}_i = 0$ with probability $1-x_i$. Define the set function:
$$\displaystyle f_{\mbs{y}}(S) = \max_{T \subseteq S, |T| \leq k}\sum_{i \in T}c_iy_i, \ \forall S \subseteq N.$$ For a fixed $\mb{y} \in {\cal Y} \subseteq \{0,1\}^n$, the function $f_{\mbs{y}}:2^{N} \rightarrow \mathbb{R}_+$  is a monotone submodular function and the corresponding value of the concave closure $f^+_{\mbs{y}}(\mb{x})$ is computable in polynomial time \cite{Calinescu1,natbook}.
The distributionally robust $k$-sum combinatorial optimization problem is then formulated as:
$$\min_{\mbs{y} \in {\cal Y} \subseteq \{0,1\}^n}f^+_{\mbs{y}}(\mb{x}),$$
where the minimization is over the worst-case dependence of the random weights. The optimal solution for this problem (denoted by $\mb{y}^+$) can be found in polynomial time under the assumption that optimizing a linear function over the set ${\cal Y}$ is possible in polynomial time (see Chapter 2 in \cite{natbook}). If the random weights are mutually independent, the stochastic optimization problem is formulated as: $$\min_{\mbs{y} \in {\cal Y} \subseteq \{0,1\}^n}F_{\mbs{y}}(\mb{x}).$$ Using Theorem \ref{thm:calinescushipra}, we obtain a $e/(e-1)$ approximation algorithm for the stochastic optimization problem since:
$$\displaystyle F_{\mbs{y}^+}(\mb{x})  \leq f_{\mbs{y}^+}^+(\mb{x}) = \min_{\mbs{y} \in {\cal Y} \subseteq \{0,1\}^n}f_{\mbs{y}}^+(\mb{x})\leq (e/(e-1))\min_{\mbs{y} \in {\cal Y} \subseteq \{0,1\}^n}F_{\mbs{y}}(\mb{x}).$$ Now suppose the random weights are pairwise independent. The distributionally robust $k$-sum combinatorial optimization problem given a fixed marginal probability vector $\mb{x}$ and assuming pairwise independence is formulated as:
$$\min_{\mbs{y} \in {\cal Y} \subseteq \{0,1\}^n}f^{++}_{\mbs{y}}(\mb{x}).$$
The results derived in this paper imply a $4/3$ approximation algorithm for this problem in the following two cases: (a) for $n = 3$ and general values of $\mb{x}$, and (b) for general $n$ with all small or all large values in $\mb{x}$. Specifically: $$\min_{\mbs{y} \in {\cal Y} \subseteq \{0,1\}^n}f^{++}_{\mbs{y}}(\mb{x}) \leq f^{++}_{\mbs{y}^+}(\mb{x})  \leq 4/3\min_{\mbs{y} \in {\cal Y} \subseteq \{0,1\}^n}f^{++}_{\mbs{y}}(\mb{x}).$$
\end{example}
\begin{example}
\textit{Properties of the continuous relaxation of the weighted coverage function}: Let $w_1,\ldots,w_m$ denote nonnegative weights and $T_1,\ldots,T_n$ be subsets of $\{1,\ldots,m\}$. The weighted coverage function $f:  2^{N} \rightarrow \mathbb{R}_+$ is defined by $f(S) = \sum_{j \in \cup_{i \in S} T_i} w_j$. For each $j \in \{1,\ldots,m\}$, let $U_j = \{i \in N \ | \ j \in T_i\}$ represent the subsets that cover $j$. The multilinear extension for any $\mb{x} \in [0,1]^n$ is given by:
       \begin{equation*}
\begin{array}{lll}
\displaystyle  F(\mb{x}) & = & \sum_{j=1}^m w_j\mathbb{P}(j \mbox{ is covered by some }T_i) \\
    & = & \sum_{j=1}^m w_j(1-\prod{i \in U_j}{}(1-x_i)).
    \end{array}
    \end{equation*}
     An upper bound on $F(\mb{x})$ is obtained by using the concave closure for each individual term:
    $$\displaystyle F(\mb{x}) \leq \sum_{j=1}^m w_jf_j^+(\mb{x}) \mbox{ where } f_j^{+}(\mb{x})= \min(1,\sum_{i \in U_j}x_i).$$ The upper bound can maximized efficiently over a polytope $ {\cal P} \subseteq [0,1]^n$ in polynomial time using linear optimization or subgradient methods \cite{Calinescu1,Karimi} to obtain an optimal solution $\mb{x}^+$. Theorem \ref{thm:calinescushipra} ensures $F(\mb{x}^+)  \geq (1-1/e)\max_{\mbs{x} \in {\cal P} \subseteq [0,1]^n}F(\mb{x})$. In prior work, pipage rounding has been applied to $\mb{x}^+$ to find good solutions to the deterministic combinatorial optimization problem while preserving the quality of the approximation in expectation \cite{Karimi}. Our results imply that the solution $\mb{x}^+$ possesses an additional performance guarantee, that is, if we use the upper pairwise independent extension for each individual term then $\sum_{j=1}^{m}w_jf_j^{++}(\mb{x}^+)  \geq 3/4\sum_{j=1}^{m}w_j\max_{\mbs{x} \in {\cal P} \subseteq [0,1]^n}f_j^{++}(\mb{x})$ for the cases considered in this paper.
\end{example}
We end the paper with a conjecture and invite readers to prove or disprove the conjecture.
\begin{conj}
For any $n$, any nonnegative monotone submodular function $f: 2^{N} \rightarrow \mathbb{R}_+$ and any $\mb{x} \in [0,1]^n$: $$\frac{f^{+}(\mb{x})}{f^{++}(\mb{x})} \leq \frac{4}{3}.$$
\end{conj}

\subsection*{Acknowledgement}
The research of the second author was partially supported by the MOE Tier 2 grant T2EP20124-0013 on Distributionally robust submodular optimization: Theory, algorithms and applications. We would like to thank the Associate Editor and an anonymous referee for their inputs.

\end{document}